\documentclass[graybox]{svmult}

\usepackage{type1cm}        
 \usepackage{xcolor}                           
\usepackage{makeidx}         
\usepackage{graphicx}        
\usepackage{multicol}        
\usepackage[bottom]{footmisc}

\usepackage{newtxtext}       %
\usepackage[varvw]{newtxmath}       

\makeindex             

\begin{document}

\title*{Dirichlet-Neumann Waveform Relaxation Method with Multiple Subdomains for Reaction-Diffusion Equation with a Time Delay \thanks{Accepted for proceedings in the 29th International Conference on Domain Decomposition Methods (DD29), 2026.}}
\author{Bankim C. Mandal\orcidID{0009-0009-0134-0422} and\\ Deeksha Tomer\orcidID{0009-0002-3269-883X}}
\titlerunning{DNWR with Multiple Subdomains for Reaction-Diffusion Equation with Time Delay}
\institute{Bankim C. Mandal \at IIT Bhubaneswar, \email{bmandal@iitbbs.ac.in}
\and Deeksha Tomer \at IIT Bhubaneswar \email{a21ma09002@iitbbs.ac.in}}

%
%
\maketitle



 \abstract*{In this study, we present the numerical investigation of the Dirichlet–Neumann Waveform Relaxation (DNWR) algorithm applied to multiple subdomains for the reaction–diffusion equation with time delay. Various arrangements of transmission conditions between subdomains are explored and a series of numerical experiments are conducted to evaluate and compare the efficiency and effectiveness of these configurations.}

\section{Introduction}
\label{sec:1}

The reaction–diffusion equation with time delay can be used to model a wide range of real-world phenomena, including predator–prey interactions, population dynamics, and the spatial spread of bacterial or viral diseases. In case of disease transmission, the delay term represents incubation or latency periods.

In neuroscience, it models neural activity in the central nervous system, where information exchange between neurons involves finite transmission delays \cite{coombes2005waves}.
Continuous combustion systems  
can also be modeled using delayed reaction–diffusion frameworks \cite{Wu}.
Introducing time delays into partial differential equation models allows researchers to capture important dynamic behaviors that may otherwise be overlooked. This modeling approach offers valuable insights for better understanding, predicting, and controlling complex phenomena across various scientific and engineering fields. Given the broad scope and computational demands of such applications, developing efficient parallel solution strategies becomes essential.

  Parallel computing techniques \cite{gander201550} are crucial for solving large-scale numerical problems that frequently arise in physics and engineering. A prominent and widely adopted strategy is domain decomposition \cite{bjor,widlund,tang}, which employs a divide-and-conquer technique. In this method, the computational domain is divided into multiple subdomains, which may be overlapping or non-overlapping. This framework naturally allows for the use of different numerical methods in each subdomain, making it well suited to handle the varying characteristics present in physical phenomena, while also tackling the computational challenge of solving very large systems.
  
  In this work, we extend the Dirichlet-Neumann Waveform Relaxation (DNWR) \cite{gander2013dirichlet,sana2024convergence} algorithm to the reaction-diffusion equation with time delay over multiple subdomains. DNWR is an extension of the steady-state Dirichlet–Neumann (DN) algorithm for the evolution problem. The proposed approach generalizes the steady-state DN method to the time-dependent setting by iteratively updating interface data through a convex combination of the current Dirichlet solution and the previous Neumann trace. Earlier work \cite{mandal2025dirichlet} demonstrated two-step convergence for the delayed reaction-diffusion problem of two subdomains. The numerical investigation of the DNWR method for the multi-subdomain case has not been comprehensively addressed previously, and we present it in this work.
    
\section{DNWR for Multiple Subdomains - Numerical Experiments}

For demonstration purposes, the following reaction-diffusion equation that incorporates constant time delay (as presented in \cite{shulin,mandal2025dirichlet}) in $\Omega\subset \mathbb{R}^d$ is examined.
\begin{equation}
\begin{array}{rl}\label{model_eq}
\partial_t w-\upsilon ^2\Delta  w+a_1w(\bar{x}, t)+a_2w(\bar {x}, t-\tau )&=f(\bar x,t),\ \ \ (\bar x, t)\in \Omega \times (0, T), \\ 
 w(\bar x, t) &=w_0(\bar x, t), \ \ \  (\bar x, t) \in \bar{\Omega}\times [-\tau , 0], \\ 
 w(\bar x , t ) &=0, \ \ \ \bar x\in \partial \Omega, t\in (0, T) , 
\end{array}
\end{equation}
where  $ \upsilon>0$, $\tau>0$ and the coefficients $a_1, a_2$ are arbitrary considering $a_2\neq 0$.
DNWR algorithm is applied in multiple subdomains $\Omega_i, i=1,2 \ldots 2n+1.$ The decomposition is done in the form of strips. For ($d=1$) decompositions are in Fig. \ref{fig3_5}, \ref{arrangement_2}, and \ref{fig:my_label}, for different arrangements. Let the interface boundaries be $\Gamma_i=\partial \Omega_i\cap \partial \Omega_{i+1}$. The interface functions $h_i^0(\bar x,t)$ on $\Gamma_i \times (0,T)$ is initialized randomly. For numerical experiments, we initialize all interface functions as
$h_i^0(\bar{x}, t) =h^0_i= t^2; $ and we work on the error equation ($L^2(0,T;L^\infty(\Omega))$ norm) for the equation \eqref{model_eq} to observe error convergence on the interface boundaries. Thus, the boundary values for the error equation are $h_0^0=h^0_{2n+1}=0$. For numerical simulations (see Fig.~\ref{arrangement1}, \ref{arrangement2}, and \ref{arrangement_3}), a 1-D spatial domain $\Omega = (0, 5)$ is divided into five equal subdomains, $\Omega_i$, $i \in \{1, \ldots, 5\}$. 
Parameters are $a_1 = 0$, $a_2 = 0.028$ and delay $\tau = 3$ for longer time interval $[0, 10]$. In contrast, we use $\tau=0.03$ for experiments with a shorter time window $[0, 0.1],$ with $a_1$, $a_2$ unchanged. The numerical discretization is carried out using the implicit finite difference scheme with mesh size $\Delta x = 0.1$ and time step $\Delta t = 0.2$. 
      Depending on the positioning of transmission conditions (Dirichlet and Neumann \cite{gander2013dirichlet}) at the interface boundaries, 3 possible arrangements arise in the multi-subdomain case for DNWR. In this section, we investigate these arrangements numerically.
\subsection{Arrangement 1}
\label{subsec:1}
\begin{figure}[!ht]
     \sidecaption
     \includegraphics[width=7cm]{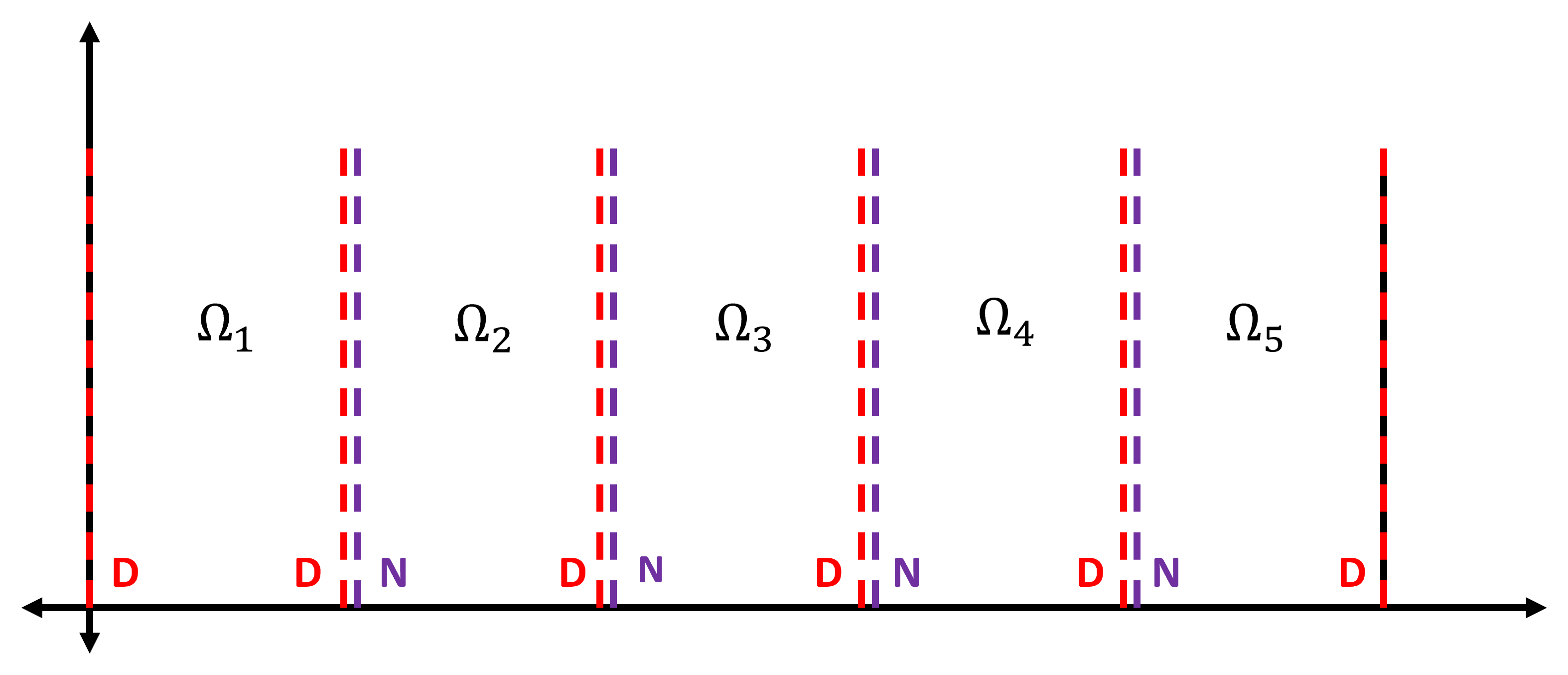}
     \caption{Arrangement 1 of Boundary conditions for DNWR in multisubdomain setup}
     \label{fig3_5}
 \end{figure}
\begin{figure}[!ht]
    \centering
    {\includegraphics[width=0.462\linewidth]{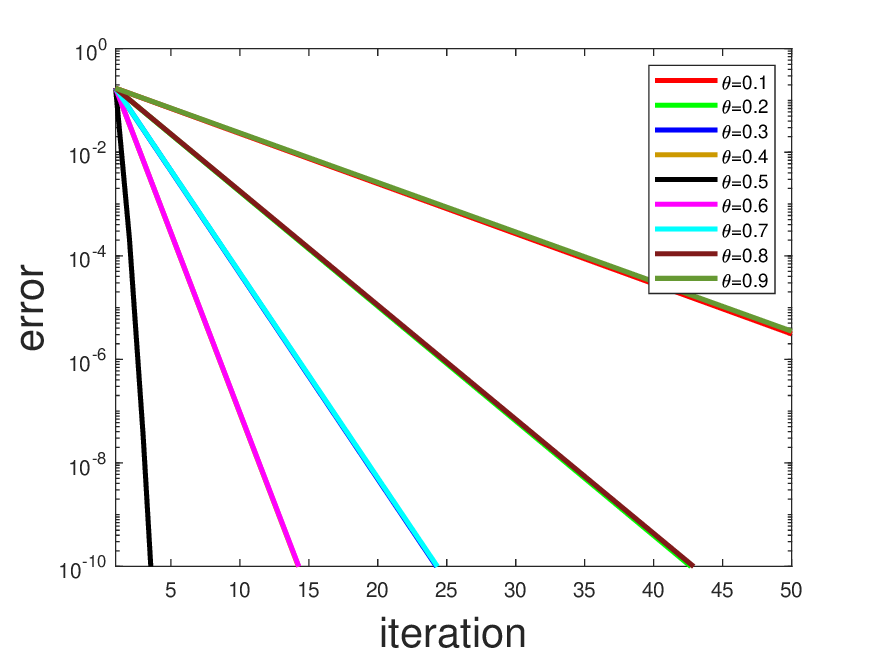} }%
    \qquad
    {\includegraphics[width=0.462\linewidth]{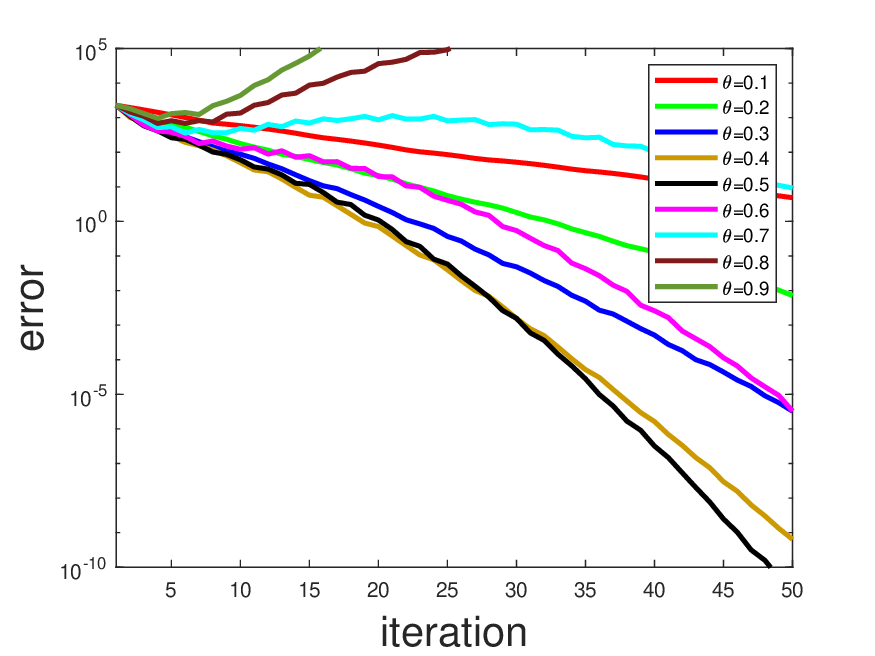} }%
    \caption{Convergence result of DNWR algorithm for arrangement $1$. Left: DNWR for shorter time window; Right: DNWR for longer time window}
    \label{arrangement1}
    \end{figure}

The formulation for two subdomains \cite{mandal2025dirichlet} is extended to multiple subdomains, as illustrated in Fig.~\ref{fig3_5}.
Starting with an initial guess along the boundary, a Dirichlet subproblem is solved within the first subdomain. Subsequently, the subproblems in the remaining subdomains are solved using mixed Neumann–Dirichlet boundary conditions and the numerical results are reported in Fig.~\ref{arrangement1}. The corresponding DNWR algorithm on the error equations for the reaction-diffusion equation with time delay in Arrangement~1 for $k=1,2,\ldots$ are given as follows:
    
     \begin{equation}
\left\{\begin{array}{rl}
\partial_t e_1^k-\upsilon ^2\Delta e_1^k+a_1e_1^k(\bar x, t)+a_2e_1^k(\bar x, t-\tau ) &=0, \ \ \ (\bar x, t)\in \Omega_1\times (0, T),  \\ 
  e_1^k(\bar x, t)&=0, \  \ \ (\bar x, t)\in \bar\Omega_1\times (-\tau, 0), \\ 
e_1^k&=0, \ \ \ \text{on} \ \ \partial \Omega \cap \partial \Omega_1\times (0,T),\\ 
e_1^k&=h_1^{k-1}, \ \ \  \text{on} \ \Gamma_1 \times (0,T).
\end{array}\right.
\end{equation}

For $j=2,3 \ldots 2n+1$
   \begin{equation}
\left\{\begin{array}{rl}
\partial_t e_j^k-\upsilon ^2\Delta e_j^k+a_1e_j^k(\bar x, t)+a_2e_j^k(\bar x, t-\tau )&=0, \ (\bar x, t)\in \Omega_j\times (0, T),  \\ 
  e_j^k(\bar x, t)&=0, \  (\bar x, t)\in \bar\Omega_j\times (-\tau, 0), \\ 
  \partial_{\mathbf{n}_{j,j-1}} e_j^k&=-\partial_{\mathbf{n}_{j-1,j}} e_{j-1}^k, \ \ \text{on} \ \Gamma_{j-1}\times(0,T),\\
e_j^k&=h_j^{k-1},\ \ \ \ \ \ \   \ \ \ \  \ \text{on} \ \ \Gamma_j\times(0,T).
\end{array}\right.
\end{equation}  
The interface update condition, with $\theta$ denoting the relaxation parameter is as follows (for all three arrangements,  $\theta \in (0,1]$):

\begin{equation}
h_j^{k}(\bar x,t)=\theta e_{j+1}^{k}|_{\mathbf{\Gamma}_{j}\times(0,T)} +(1-\theta) h_j^{k-1}(\bar x, t).
\end{equation}

    \subsection{Arrangement $2$}
        \begin{figure}[!ht]
     \sidecaption
     \includegraphics[width=7 cm]{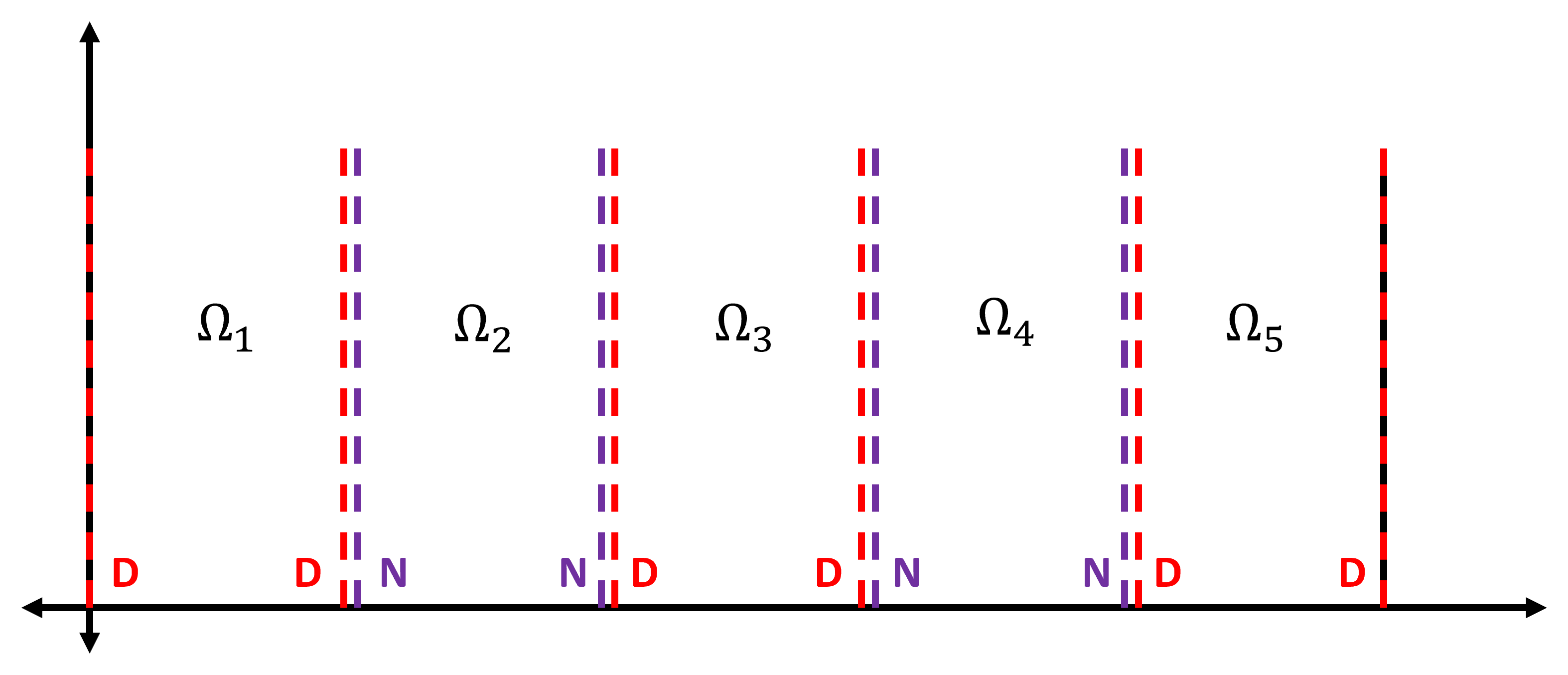}
     \caption{Arrangement 2 of boundary conditions for DNWR in multisubdomain setup}
     \label{arrangement_2}
 \end{figure}
 
\begin{figure}[!ht]
    \centering
    \includegraphics[width=0.465\linewidth]{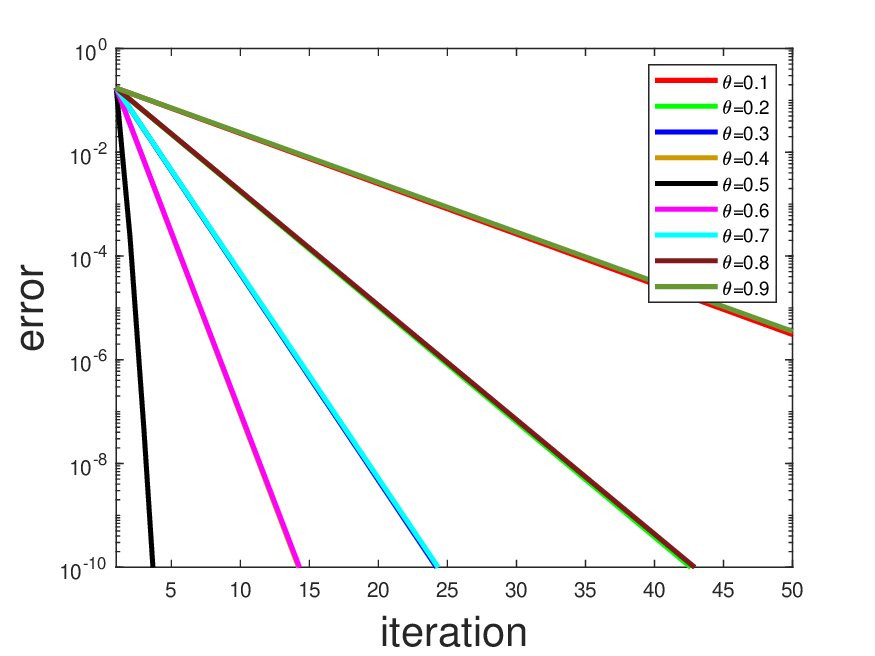} %
    \qquad
    \includegraphics[width=0.465\linewidth]{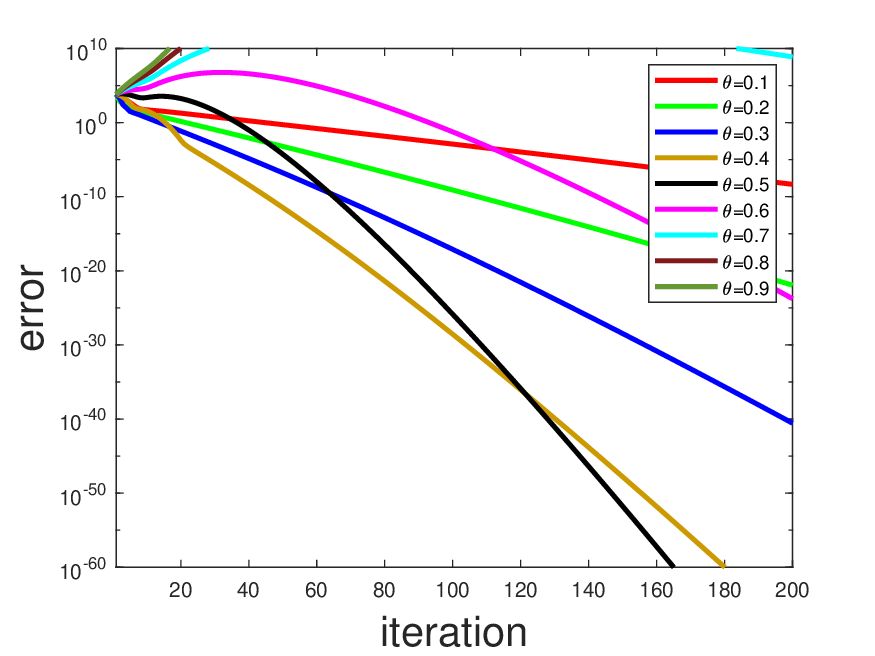} %
    \caption{ Convergence result of DNWR algorithm for arrangement $2$. Left: DNWR for shorter time window; Right: DNWR for longer time window.}
    \label{arrangement2}
    \end{figure}

    We first solve the Dirichlet subproblems in parallel over the alternating domains $\Omega_1, \Omega_3, \ldots \Omega_{2n+1}$, followed by parallel solutions of the Neumann subproblems in the remaining domains $\Omega_2, \Omega_4 \ldots \Omega_{2n}$ (see Fig.~\ref{arrangement_2} and for numerical result see Fig.~\ref{arrangement2}). For $k=1,2,\ldots$ and $j = 1, 3, \ldots 2n+1$,
    
\begin{equation}
\left\{\begin{array}{rl}
\partial_t e_j^k-\upsilon ^2\Delta e_j^k+a_1e_j^k(\bar x, t)+a_2e_j^k(\bar x, t-\tau )&=0, \ (\bar x, t)\in \Omega_j\times (0, T),  \\ 
  e_j^k(\bar x, t)&=0, \  (\bar x, t)\in \bar\Omega_j\times (-\tau, 0), \\ 
e_j^k&=h^{k-1}_{j-1},  \  \text{on} \ \Gamma_{j-1}\times(0,T),\\ 
e_j^k&=h_j^{k-1},  \  \text{on} \ \Gamma_j\times(0,T),
\end{array}\right.
\end{equation}  
 for $i=2,4 \ldots 2n$,
   \begin{equation}
\left\{\begin{array}{rl}
\partial_t e_i^k-\upsilon ^2\Delta e_i^k+a_1e_i^k(\bar x, t)+a_2e_i^k(\bar x, t-\tau )&=0, \ (\bar x, t)\in \Omega_i\times (0, T),  \\ 
  e_i^k(\bar x, t)&=0, \  (\bar x, t)\in \bar\Omega_i\times (-\tau, 0), \\
  \partial_{\mathbf{n}_{i,i-1}} e_i^k&=-\partial_{\mathbf{n}_{i-1,i}} e_{i-1}^k, \ \ \ \  \text{on} \ \Gamma_{i-1}\times(0,T),\\
  \partial_{\mathbf{n}_{i,i+1}} e_i^k&=-\partial_{\mathbf{n}_{i+1,i}} e_{i+1}^k, \ \ \ \  \text{on} \ \Gamma_{i}\times (0,T).\\
\end{array}\right.
\end{equation}
At the interface, the update rule involving the relaxation parameter $\theta$ is defined as:
\begin{equation}
h_j^{k}(\bar x,t)=\theta e_{j+1}^{k}|_{\mathbf{\Gamma}_{j}\times(0,T)} +(1-\theta) h_j^{k-1}(\bar x, t);\ j=1,3\ldots 2n+1,
\end{equation}
\begin{equation}
h_i^{k}(\bar x,t)=\theta e_{i}^{k}|_{\mathbf{\Gamma}_{i}\times(0,T)} +(1-\theta) h_i^{k-1}( \bar x,t);\ \ i=2,4 \ldots 2n.
\end{equation}

    \subsection{Arrangement $3$}

        \begin{figure}[!ht]
     \sidecaption
     \includegraphics[width=7cm]{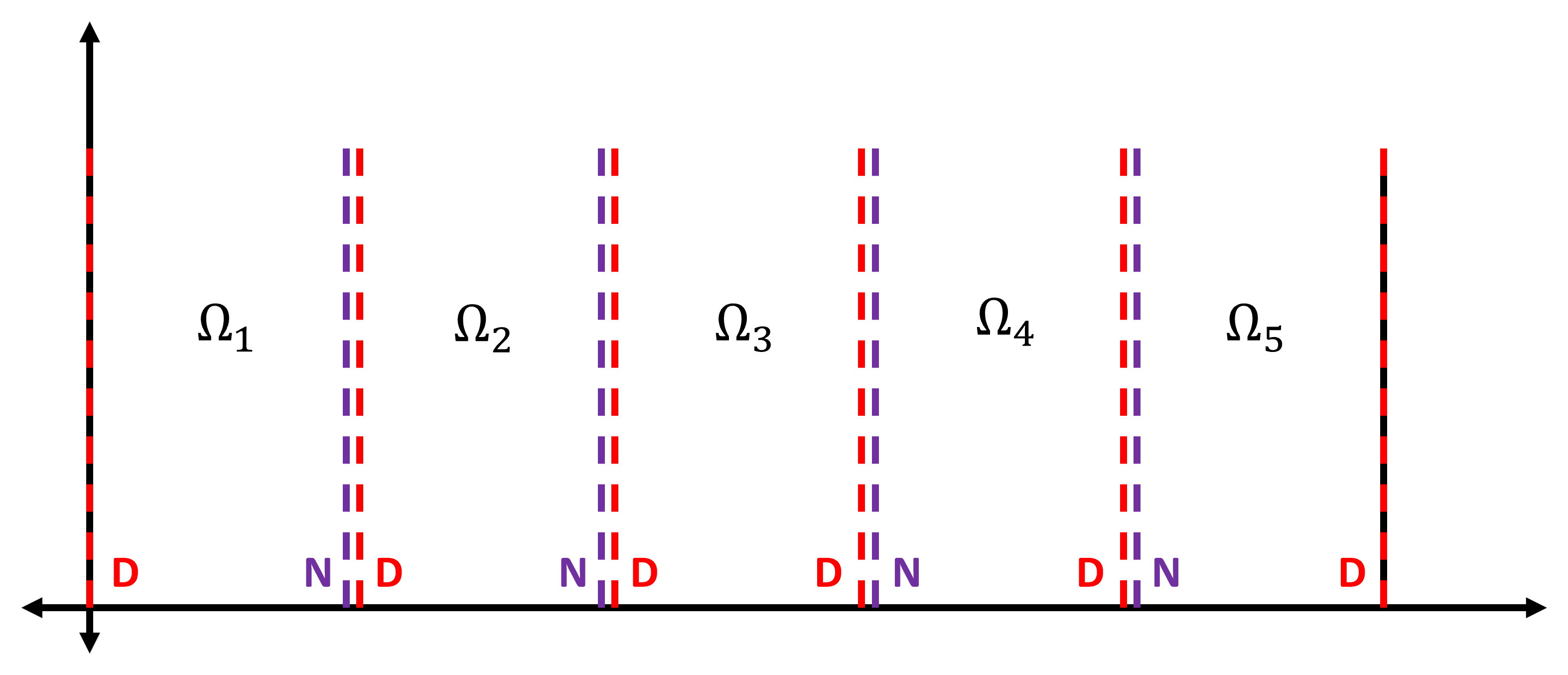}
     \caption{Arrangement 3 of boundary conditions in multisubdomain setup}
     \label{fig:my_label}
 \end{figure}
\begin{figure}[!ht]
    \centering
    {\includegraphics[width=0.462\linewidth]{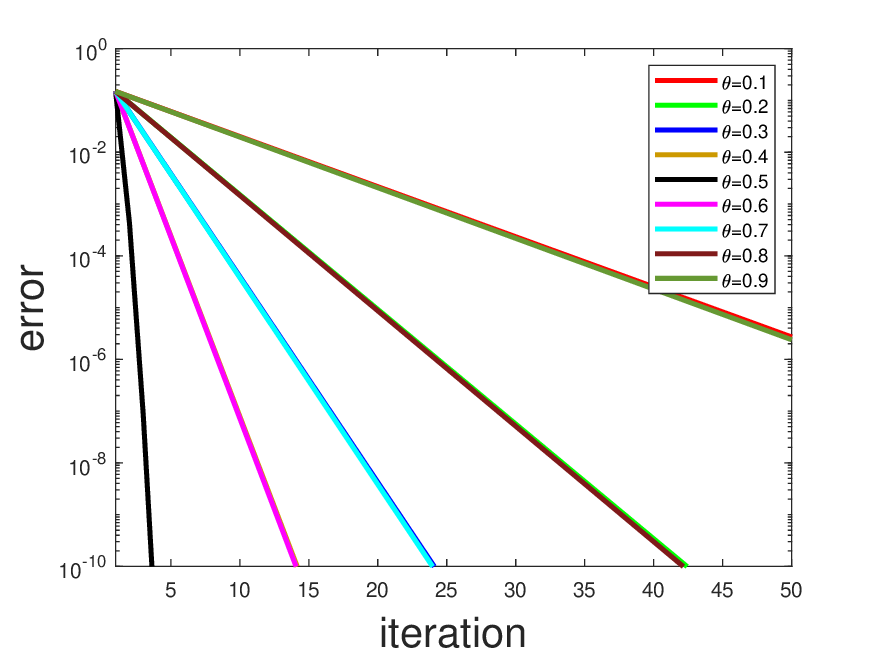} }%
    \qquad
    {\includegraphics[width=0.462\linewidth]{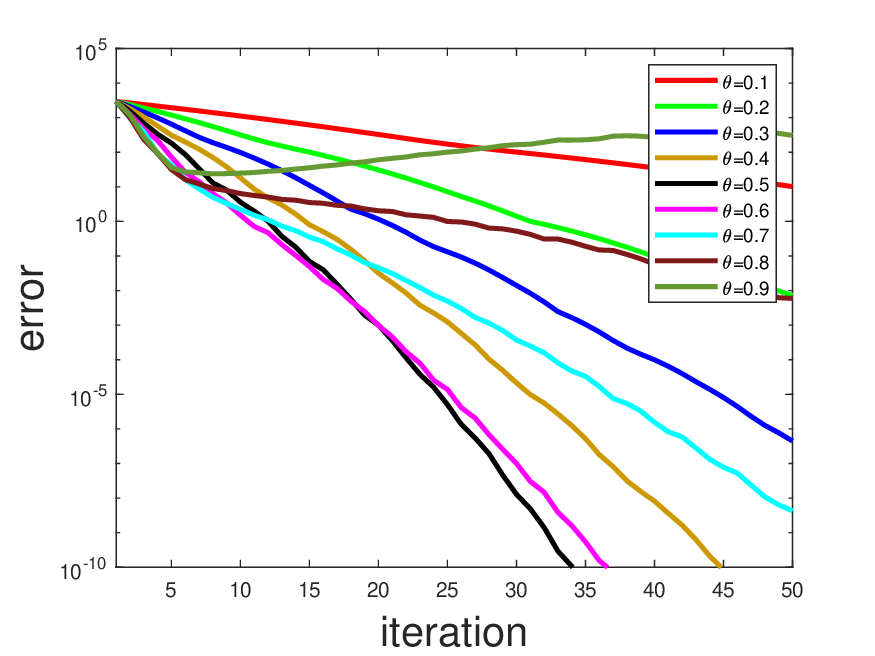} }%
    \caption{ Convergence result of DNWR algorithm for arrangement $3$. Left: DNWR for shorter time window; Right: DNWR for longer time window. }
    \label{arrangement_3}
    \end{figure}

    In this arrangement, the Dirichlet problem is first solved in the central subdomain, followed by solving the Dirichlet-Neumann problem in the neighboring subdomains subsequently (see Fig. \ref{fig:my_label}). \textcolor {black}{The numerical result is in Fig. \ref{arrangement_3}}.   
    For $k=1,2,\ldots$ 

\begin{equation}
\resizebox{0.95\linewidth}{!}{%
    $ 
    \left\{\begin{array}{rl}
    \partial_t e_{n+1}^k-\upsilon ^2\Delta  e_{n+1}^k+a_1e_{n+1}^k(\bar{x}, t)+a_2e_{n+1}^k(\bar {x}, t-\tau )&=0, \ (\bar x, t)\in \Omega_{n+1}\times (0, T),  \\ 
     e_{n+1}^k(\bar x, t)&=0, \  (\bar x, t)\in \bar\Omega_{n+1}\times (-\tau, 0), \\ 
    e_{n+1}^k&=0, \  \text{ on}\ \partial\Omega \cap\partial\Omega_{n+1} \times (0,T),\\ 
    e_{n+1}^k&=h_l^{k-1},  \ \text{on}\ \Gamma_l \times (0,T), l=n,n+1.
    \end{array}\right.
    $ 
}
\end{equation}
For $n \geq i\geq 1$, in this reverse ordering, we solve 
\begin{equation}
\left\{\begin{array}{rl}
\partial_t e_i^k-\upsilon ^2\Delta e_i^k+a_1e_i^k(\bar x, t)+a_2e_i^k(\bar x, t-\tau )&=0, \ \ \ \ \ \ \  (\bar x, t)\in \Omega_i\times (0, T),  \\ 
  e_i^k(\bar x, t)&=0,  \ \ \ \ \ \ \ \ (\bar x, t)\in \bar\Omega_i\times (-\tau, 0), \\ 
e_i^k&= h^{k-1}_{i-1}, \ \  \ \ \ \ \ \ \ \ \text{on} \  \Gamma_{i-1}  \times (0,T),\\
\partial_{\mathbf{n}_{i,i+1}} e_i^k&=-\partial_{\mathbf{n}_{i+1,i}} e_{i+1}^k, \ \ \ \text{on} \ \Gamma_i \times (0,T).
\end{array}\right.
\end{equation}
Finally, for $n+2\leq j \leq 2n+1$
\begin{equation}
\left\{\begin{array}{rl}
\partial_t e_j^k-\upsilon ^2\Delta e_j^k+a_1e_j^k(\bar x, t)+a_2e_j^k(\bar x, t-\tau )&=0, \ (\bar x, t)\in \Omega_j\times (0, T),  \\ 
  e_j^k(\bar x, t)&=0, \ \ \ \ \ \ \ (\bar x, t)\in \bar\Omega_j\times (-\tau, 0), \\ 
\partial_{\mathbf{n}_{j,j-1}} e_j^k&=-\partial_{\mathbf{n}_{j-1,j}} e_{j-1}^k, \ \ \ \  \text{on} \ \Gamma_{j-1} \times (0,T),\\
e_j^k&= h^{k-1}_j,\  \ \ \ \ \ \ \ \ \ \ \  \ \ \ \text {on} \ \ \  \Gamma_{j} \times (0,T).
\end{array}\right.
\end{equation}\\
The interface update condition, with $\theta$ the relaxation parameter, is
\begin{equation}
h_i^{k}(\bar x,t)=\theta e_{i}^{k}|_{\mathbf{\Gamma}_{i}\times(0,T)} +(1-\theta) h_i^{k-1}(\bar x, t);\  1\leq i\leq n,\\
\end{equation}
\begin{equation}
h_j^{k}(\bar x,t)=\theta e_{j+1}^{k}|_{\mathbf{\Gamma}_{j}\times(0,T)} +(1-\theta) h_j^{k-1}(\bar x, t);\ n+1\leq j \leq 2n.\\
\end{equation}

\noindent\textbf{Remark:} Numerical Results (Fig. \ref{arrangement1},\ref{arrangement2},\ref{arrangement_3}) show that the third arrangement is the most efficient for longer time domain simulations, although none of the arrangements are completely parallel. Arrangement~1 is sequential, while arrangement~2,~3 are partially parallel. Arrangement~2 is well suited for parallel computing due to its easier parallelization.

    \section{Additional Numerical Experiments}
\begin{figure}[!ht]
    \centering
    \includegraphics[width=0.424\linewidth]{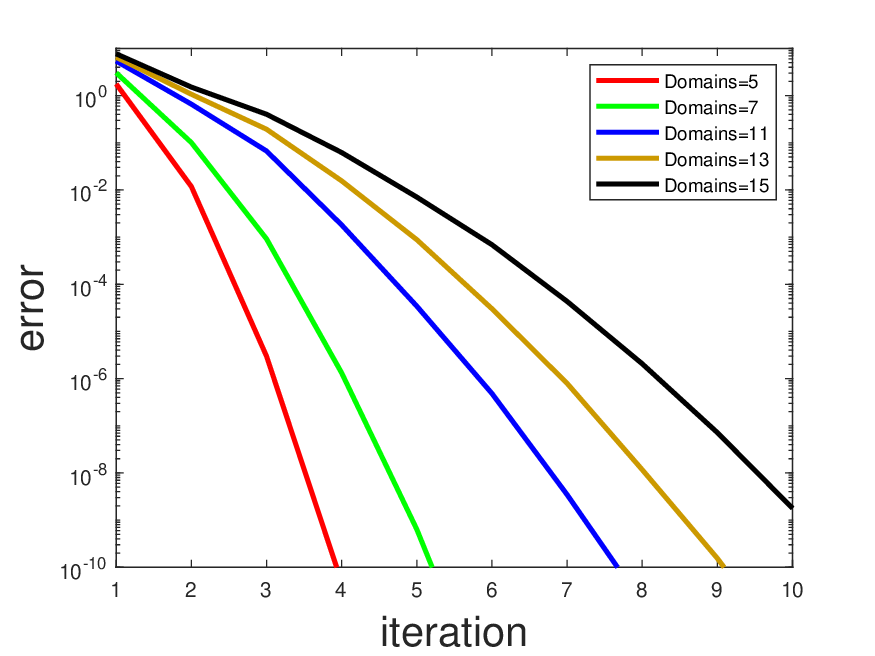} %
    \qquad
    \includegraphics[width=0.424\linewidth]{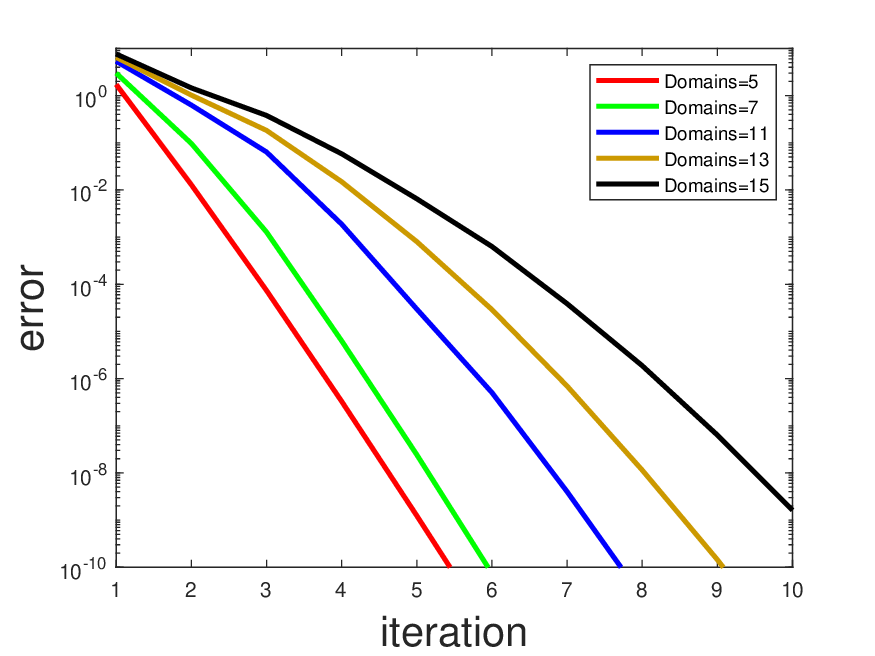} %
    \caption{DNWR convergence for different no. of subdomains when $T=0.1$, $\tau=0.03$. Left: Convergence when $a_1=0$; Right: Convergence when $a_1=1$.}
    \label{diff_no_subdom}
    \end{figure}

\begin{figure}[!ht]
    \centering
    \includegraphics[width=0.462\linewidth]{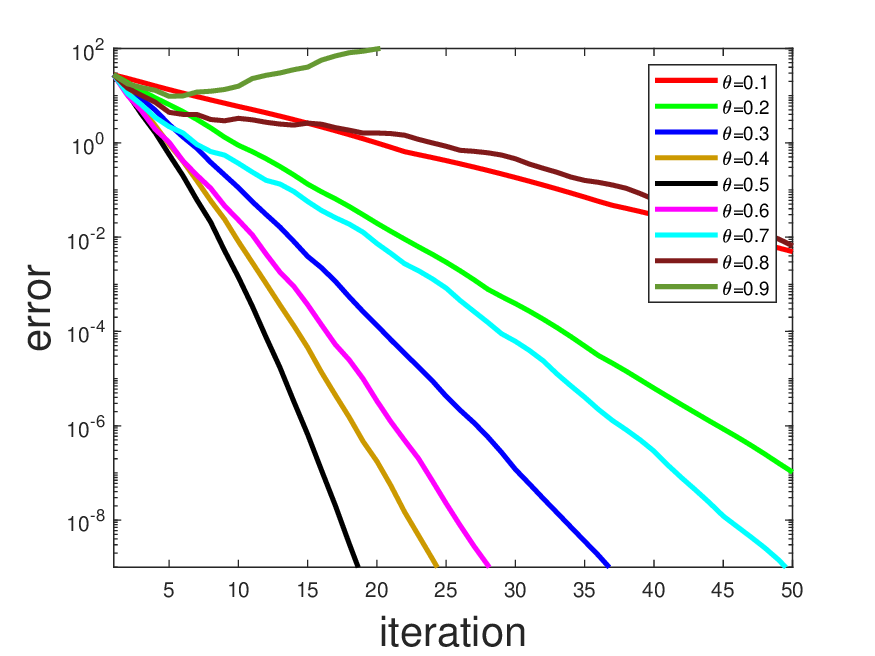} %
    \qquad
    \includegraphics[width=0.462\linewidth]{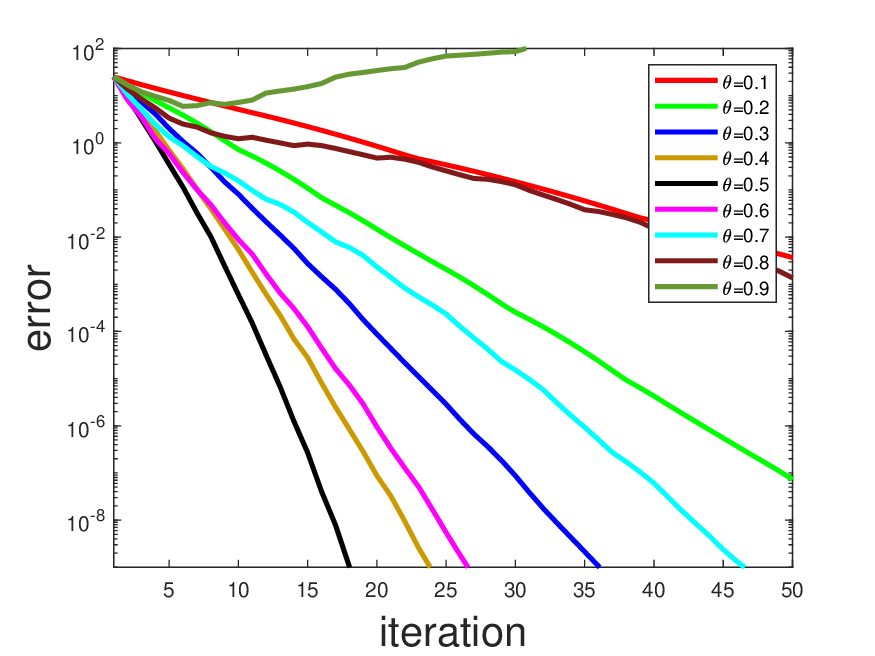} %
    \caption{ Convergence behavior of DNWR method for arrangement $3$ when $T=1$ and $\tau=0.3$ for subdomains having different sizes. Left: Convergence when $a_1=0$; Right: Convergence when $a_1=1$.}
    \label{diff_subdomain_size}
    \end{figure}

\begin{figure}[!ht]
    \centering
    \includegraphics[width=0.462\linewidth]{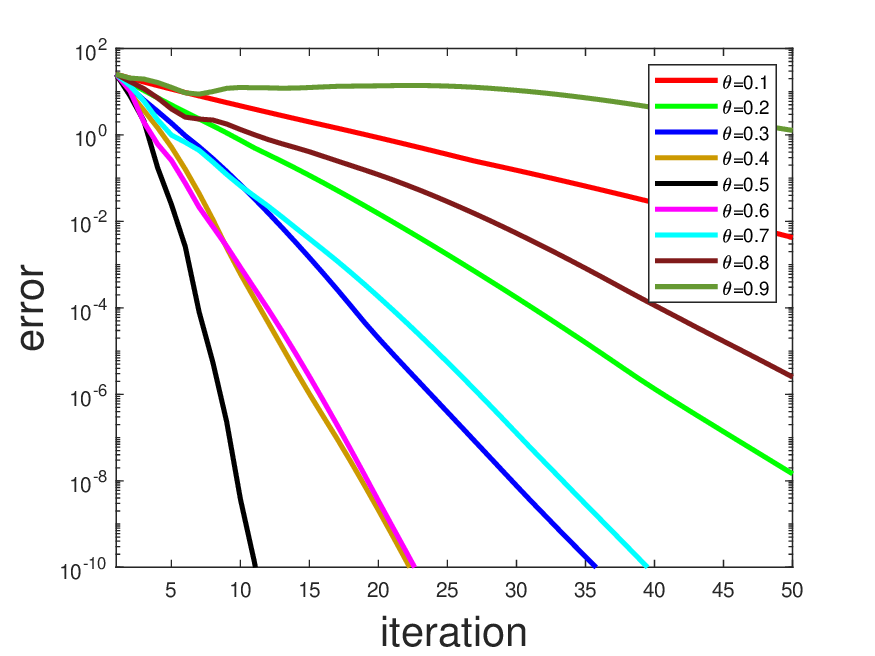}%
    \qquad 
    \includegraphics[width=0.462\linewidth]{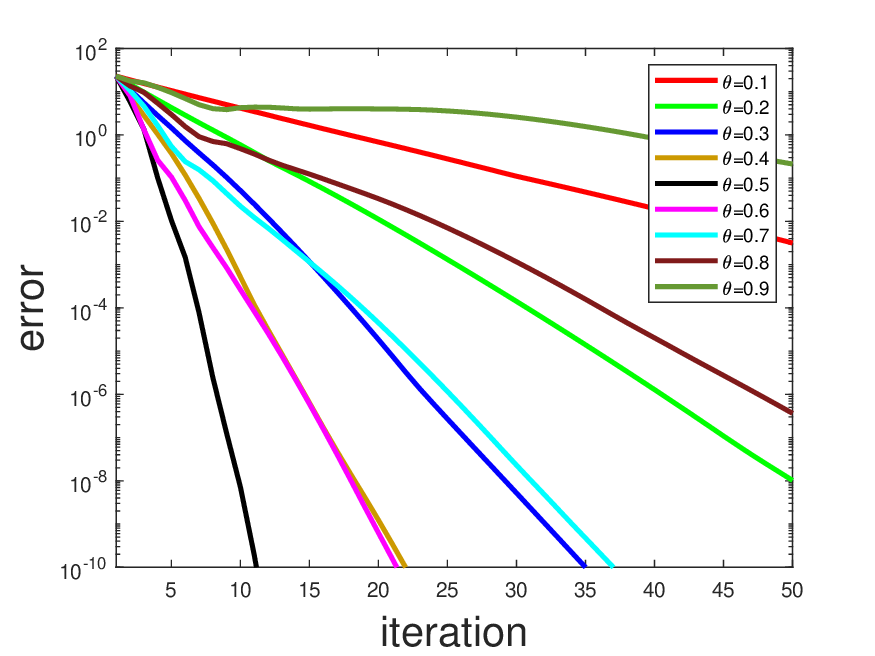}
    
    \caption{Convergence behavior of DNWR for arrangement $3$ having different interface values for $T=1$ and $\tau=0.3$. Left: Convergence when $a_1=0$; Right: Convergence when $a_1=1$.}
    \label{diff_interface}
\end{figure}

      Additional numerical experiments were performed for Arrangement~3 on error equations to further understand the behaviour of DNWR in the multisubdomain case. 
    \begin{enumerate}
    \item Numerical experiments are performed for varying numbers of subdomains with $\theta = 1/2$, $\Delta t = 0.001$, and $\Delta x$ adjusted according to the number of subdomains, as shown in Fig.~\ref{diff_no_subdom}. The results indicate that the iteration count needed to achieve convergence grows as the number of subdomains increases. 
    \item To analyze the effectiveness of Arrangement~$3$, we carried out experiments considering cases with unequal subdomain sizes, i.e., $|\Omega_1|=|\Omega_5|=1.5, |\Omega_2|=|\Omega_4|=0.5, |\Omega_3|=1 $. Fast convergence was observed for $\theta=1/2$ (Fig. \ref{diff_subdomain_size}).
    \item The performance of Arrangement~$3$ is also evaluated through experiments where different interfaces are initialized with distinct functions, namely: $(h_1^0 = t^2)$; $(h_2^0 = t)$; $(h_3^0 = \sin(t))$; and 

$h_4^0=
\{t~(0<t\le0.4),~t^2~(0.4<t\le 0.8),~\sin(t)~(0.8<t\le1)\}.$
From Fig. \ref{diff_interface} we can see that it converges for $\theta=1/2$. The performance is also examined when both subdomain sizes (as in 2) and interface values are different. Results are reported in Fig. 10, we can say that convergence is fastest for $\theta=1/2.$
 
    \end{enumerate}

     \begin{figure}[!ht]
    \centering
    \includegraphics[width=0.46\linewidth]{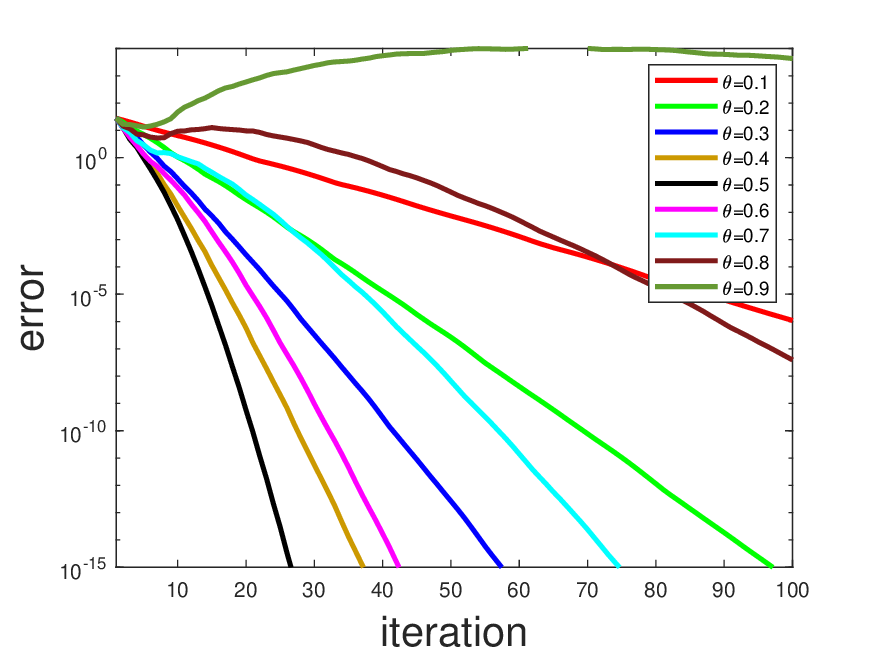} %
    \qquad
    \includegraphics[width=0.46\linewidth]{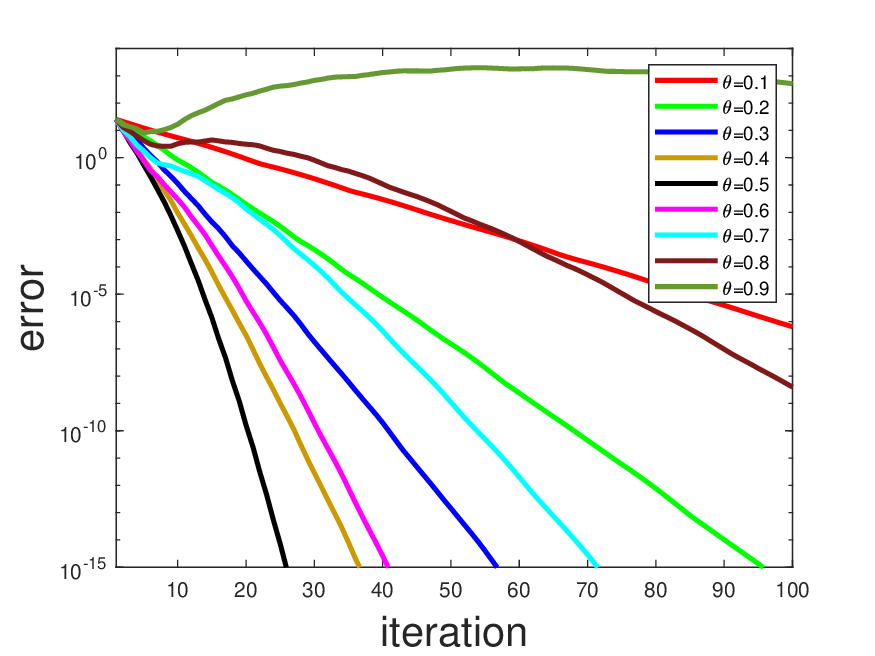} %
    \caption{ DNWR convergence behavior for arrangement $3$ when $T=1$ and $\tau=0.3$ for subdomains having different sizes and different interface values. Left: Convergence when $a_1=0$; Right: Convergence when $a_1=1$.}
    \label{diff_subdomain_size_hi}
    \end{figure}
\section{Conclusion}
\label{sec:3}

In this work, the Dirichlet–Neumann Waveform Relaxation (DNWR) algorithm has been extended to handle multiple subdomains. Results from our numerical experiments demonstrate that the third arrangement outperforms the other configurations in terms of iteration efficiency. We also obtained the fastest convergence for $\theta=1/2$ across various interface conditions and subdomain sizes. However, with the increase number of subdomains, the iteration counts required for convergence also grows.

 \bibliographystyle{spmpsci}
\bibliography{reference}

\end{document}